\documentclass[twoside,notitlepage,12pt]{article}
\usepackage{graphicx}
\usepackage{srcltx}

\usepackage{amssymb}
\usepackage[leqno]{amsmath}
\usepackage{amsfonts}
\usepackage{amsopn}
\usepackage{amstext}
\usepackage{amsthm}

\usepackage[all]{xy}
\newdir{ >}{{}*!/-9pt/@{>}}

\usepackage{verbatim}
\usepackage[colorlinks]{hyperref}

\providecommand{\cal}{\mathcal}
\renewcommand{\Bbb}{\mathbb}

\newenvironment{pf}{\begin{proof}}{\end{proof}}

\newcommand{\Ef}{{\cal{F}}}

\newcommand{\Pee}{{\cal{P}}}

\newcommand{\R}{{\Bbb{R}}}
\newcommand{\N}{{\Bbb{N}}}

\newcommand{\al}{\alpha}

\newcommand{\sig}{\sigma}
\newcommand{\eps}{\varepsilon}
\newcommand{\st}{such that}
\renewcommand{\phi}{\varphi}
\renewcommand{\rho}{\varrho}

\newcommand{\rest}{\restriction}

\newcommand{\ntr}{{n\in\omega}}

\newcommand{\loe}{\leq}
\newcommand{\goe}{\geq}

\newcommand{\subs}{\subseteq}
\newcommand{\sups}{\supseteq}

\newcommand{\arre}{\par \medskip}

\newcommand{\poset}{{\Bbb{P}}}

\newcommand{\meet}{\wedge}

\newcommand{\join}{\vee}

\newtheorem{tw}{Theorem}[section]

\newtheorem{lm}[tw]{Lemma}
\newtheorem{prop}[tw]{Proposition}
\newtheorem{claim}[tw]{Claim}
\theoremstyle{definition}
\newtheorem{df}[tw]{Definition}
\newtheorem{ex}[tw]{Example}

\theoremstyle{remark}
\newtheorem{uwgi}[tw]{Remark}

\newcommand{\setof}[2]{\{#1\colon #2\}}

\newcommand{\seq}[1]{\langle #1 \rangle}

\newcommand{\sett}[2]{\{#1\}_{#2}}
\newcommand{\sn}[1]{\{#1\}} 
\newcommand{\dn}[2]{\{#1,#2\}} 
\newcommand{\pair}[2]{\langle #1, #2 \rangle} 
\newcommand{\triple}[3]{\langle #1, #2, #3 \rangle} 
\newcommand{\map}[3]{#1\colon #2 \to #3} 

\newcommand{\dsl}[1]{{\mathbb S(#1)}} 

\newcommand{\fin}[1]{[#1]^{<\omega}}

\newcommand{\Pri}{{\mathbb P}}

\providecommand{\nat}{\omega}

\newcommand{\cee}[1]{{C}\!\left({#1}\right)}
\newcommand{\ceep}[1]{{\mathcal C_p}\!\left({#1}\right)}

\newcommand{\osc}{\operatorname{osc}}

\newcommand{\sld}{{\bf SLD}}

\title{Topological properties of the continuous function spaces on some ordered compacta}
\author{W. Kubi\'s\footnote{The first author was
supported by the GA\v{C}R grant P
201/12/0290 (Czech Republic) and by a visiting position at Department of Mathematical Analysis, University of Valencia, Spain (April -- October 2010).}, A. Molt\'o\footnote{The second author was
supported by MTM2011--22457, Ministerio de Econom\'{\i}a y Competitividad (Spain).}~ and S. Troyanski\footnote{The third author was
supported by MTM2011--22457, Ministerio de Econom\'{\i}a y Competitividad (Spain), by MCI MTM2011-22457, by Science Foundation
Ireland SFI11/RFP.1/MTH/3112, and a project of the Institute of Mathematics and Informatics,
Bulgarian Academy of Science.
}
}

\begin{document}

\maketitle
\begin{abstract}
Some new classes of compacta $K$ are considered for which $C(K)$ endowed with the pointwise topology has a countable cover by sets of small local norm--diameter.
\end{abstract}


\section{Introduction}

A topological notion introduced in \cite{JNRMat} plays an important role in the study of topological and renorming properties of Banach spaces \cite{Rajakade}: If $(X,{\mathcal T})$ is a topological space and $\rho$ is a metric on $X$, we say that it has a {\em countable cover by sets of small local $\rho$--diameter} if for every $\eps >0$ we can write $X=\bigcup_{n\in \N}X_n$ in such a way that for any $n\in\N$ and every $x\in X_n$ there exists a ${\mathcal T}$--open set $U$ \st\ $x\in U$ and $\rho$--diam$(U\cap A)<\eps$.
Here we consider some compacta $K$ \st\ $C(K)$, endowed with the pointwise topology, has a countable cover by sets of small local $\rho$--diameter when $\rho$ is the norm--metric, or $C_{p}(K)$ has \sld\ for short.  Let us recall that given a topology $\tau$, coarser than the norm topology on a Banach space $X$, we say that $X$ has $\tau$--Kadets when $\tau$ and its norm topology coincide in the unit sphere.  If $C(K)$ has a pointwise Kadets equivalent norm then $C_{p}(K)$ must have \sld\ \cite{JNRMat}, that in turn implies that $C_{p}(K,\{0,1\})$ is $\sigma$--discrete.  Whether any of the converse implications holds is a well-known open problem.  However M. Raja has shown, roughly speaking, that \sld\ is {\em very close} to the existence of a pointwise Kadets renorming, namely from \cite{Rajakade} it follows that if $X$ is a Banach space endowed with a topology $\tau$, coarser than the norm topology, then $(X,\tau )$ has the property \sld\ for the norm if, and only if, there exists a non negative symmetric homogeneous $\tau$--lower semicontinuous function (that may be not convex) $F$ on $X$ with $\| \cdot \| \leq F\leq 3\| \cdot \|$ such that the norm topology and $\tau$ coincide on the set $\{ x\in X:\ F(x)=1\}$.  In \cite{RibaSt00} it is proved that $C_{p}(K\times L)$ has \sld\ whenever $C_{p}( L)$ has it and $C_{p}(K)$ has a pointwise Kadets norm \cite{RibaSt00}.  Let us mention that if $X$ is a Banach space \st\ $(X,\hbox{\rm weak})$ has the \sld\ property for the norm and the bidual of $(X,\| \cdot \| )$ is strictly convex then $X$ has a locally uniformly rotund equivalent norm \cite{MOTPL}.  (A norm
$\|\cdot\|$ in a Banach space is locally uniformly rotund if $\lim_{k}\left\| x_{k}-x\right\| =0$ whenever $\lim_{k} \left\| \left(x_{k}+x\right) /2 \right\| =  \lim_{k}  \left\|x_{k}
\right\| =\|x\|$.)  Despite no topological characterization has been obtained for those $K$'s \st\ $C_{p}(K)$ has \sld , some light on this questions has been shed in some particular classes of compacta \cite{Haydon}, \cite{HJNR}, \cite{BurKuTod}, \cite{HMO}.

In this note we present two classes of compact spaces $K$ for which the spaces $C_p(K)$ have \sld.
It is well-known that every compact space is a continuous image of a 0-dimensional compact space.
In turn, a 0-dimensional space can be regarded as a subspace of a Cantor cube $2^S$ which can be identified with the power-set of a fixed set $S$.
Consequently, a 0-dimensional compact space carries a partial ordering, which is just the inclusion relation.
It is natural to ask that the partially ordered compact space has the property that for every two elements $x$ and $y$ there exists their infimum $\inf \dn x y$.
Moreover, it is natural to expect that the operation $\pair x y \mapsto \inf \dn x y$ is continuous.
Once this happens, we speak about \emph{compact semilattices}.

One should not expect positive topological properties of $C(K)$ spaces, where $K$ is an arbitrary compact semilattice, since this class contains 1-point compactifications of trees studied in \cite{Haydon}.
We prove, however, that for a fairly large class of compact semilattices $K$ the space $C_p(K)$ is \sld.


\section{Preliminaries}

\newcommand{\pfilter}[1]{[#1,\rightarrow)}

A \emph{semilattice} is a partially ordered set $\pair S \loe$ which contains the minimal element (always denoted by 0) and in which every pair of elements $x,y$ has the greatest lower bound, denoted by $x \meet y$.
The element $x \meet y$ is sometimes called the \emph{meet of} $x,y$ and $S$ is sometimes called a \emph{meet semilattice}.
Some authors do not require the existence of the minimal element, we do it since we are going to consider compact semilattices in which the minimal element always exists.
A semilattice $S$ is \emph{topological} if it carries a Hausdorff topology with respect to which $\meet$ is continuous.

A \emph{filter} in a semilattice $S$ is a subset $F \subs S$ (possibly empty) satisfying $$\setof{x \in S}{a \meet b \loe x} \subs F$$ for every $a,b \in F$.
A filter $F$ is \emph{principal} if it is of the form
$$\pfilter p = \setof{x \in S}{p \loe x}.$$

Later on, we shall use some standard (although not trivial) properties of compact 0-dimensional semilattices.
For details we refer to one of the books \cite{HMS} or \cite{Gierzetal}.

In particular, we shall need the following algebraic notion.
An element $p$ of a semilattice $S$ is \emph{compact} if for every $A \subs S$ with $\sup A = p$ there exists a finite $A_0 \subs A$ such that $\sup A_0 = p$.
In particular, $0$ is a compact element.

The following fact will be used later without explicit reference:

\begin{prop}\label{Pnonexplicitsmlts}
Let $K$ be a compact 0-dimensional semilattice.
Then
\begin{enumerate}
	\item[(1)] A principal filter $\pfilter p$ is a clopen set if and only if $p$ is a compact element.
	\item[(2)] Given $a, b \in K$ such that $a \not \loe b$, there exists a compact element $p$ such that $p \loe a$ and $p \not \loe b$.
	\item[(3)] Clopen principal filters and their complements generate the topology of $K$.
	\item[(4)] Given a nonempty clopen set $A \subs K$, every minimal element of $A$ is compact.
\end{enumerate}
\end{prop}

We now make few comments concerning a Stone-like duality for semilattices.
Namely, given a compact semilattice $K$, denote by $\dsl K$ the family of all clopen filters in $K$.
By Proposition~\ref{Pnonexplicitsmlts}(1), every nonempty element $F$ of $\dsl K$ can be identified with its vertex $p$, which is a compact element such that $F = \pfilter p$.
Recall that the empty set is a filter.
Thus, $\dsl K$ can be identified with the set of all positive compact elements of $K$ plus the ``artificial" element $\infty$.
Observe that $\dsl K$, treated as the family of filters, is a semilattice in which the meet of $F, G \in \dsl K$ is $F \cap G$.
If $F = \pfilter p$, $G = \pfilter q$, then either $F \cap G = \pfilter r$, where $r = \sup \dn p q$ or $F \cap G = \emptyset$.

It turns out that this operation is reversible, namely, if $\triple S \meet 0$ is a semilattice (considered without any topology) then one can define $K(S)$ to be the family of all filters in $S$ endowed with inclusion.
The space $K(S)$ is compact 0-dimensional when endowed with the topology inherited from the Cantor cube $\Pee(S)$, the power-set of $S$.
The duality (proved in \cite{HMS}) says that $K(S)$ is canonically isomorphic to $K$ whenever $S =\dsl K$.
More precisely, given $x \in K$, define $\hat x = \setof{p\in \dsl K}{x \in p}$.
Then $\hat x$ is a filter in $\dsl K$; in other words, $\hat x \in K(\dsl K)$.
It turns out that all elements of $K(\dsl K)$ are of this form.

The second part of this note is devoted to compact distributive lattices.
Suppose that $K$ is a compact semilattice with the unique maximal element $1$.
Then, by compactness, $K$ is a complete lattice, that is, for every $A \subs K$ the $\sup A$ and $\inf A$ exist.
In fact, $\inf A$ is the limit of the net $\sett{\inf S}{S\in \fin A}$, where $\fin A$ denotes the family of all finite subsets of $A$.
On the other hand, $\sup A$ is the infimum of the set of all upper bounds of $A$ (this set contains $1$, therefore is nonempty).
We shall denote $\sup \dn x y$ by $x \join y$ (sometimes it is called the \emph{join of} $x$ and $y$).
It is natural to ask when the operation $\join$ is continuous.
Once it happens, we say that $K$ is a \emph{compact lattice}.
A lattice $\triple K \meet \join$ is \emph{distributive} if it satisfies $(a\join b)\meet c = (a\meet c) \join (b\meet c)$ for every $a,b,c \in K$.
The notion of a filter in a lattice is the same as in a semilattice.
Note that a lattice with the reversed ordering is again a lattice (where the meet is exchanged with the join).
Let us call it the \emph{reversed lattice}.
A filter in the reversed lattice will be called an \emph{ideal}.
Important for us is the notion of a prime filter.
Namely, a filter $F$ is called \emph{prime} if it is nonempty and its complement is a nonempty ideal.
As we are interested in compact 0-dimensional distributive lattices $K$, we shall work with the family of all clopen prime filters in $K$, denoted by $\poset(K)$.
This is justified by the following:

\begin{prop}
Let $K$ be a compact 0-dimensional distributive lattice.
Then for every $a, b \in K$ with $a \not \loe b$ there exists a clopen prime filter $P\subs K$ such that $a \in P$ and $b \notin P$.
\end{prop}

Fix a distributive lattice $L$.
A subset $G \subs L$ is called \emph{convex} if it is of the form $G = I \cap F$, where $I$ is an ideal and $F$ is a filter.
For convenience, we allow here that $I = L$ or $F = L$, therefore every ideal and every filter are convex sets.
Given $a,b \in L$ we define the interval $[a,b] = \setof{x \in L}{a \meet b \loe x \loe a \join b}$.
This is the minimal convex set containing $a,b$.
Assume now that $L$ is a compact distributive lattice.
Given two disjoint closed convex sets $A$, $B$, given $a_0 \in A$, there always exist $a_1 \in A$, $b_1 \in B$ such that $a_1 \in [a_0,b_1]$ and $a_1 \in [x,b_1]$. $b_1 \in [a_1,y]$ holds for every $x \in A$, $y \in B$.
Furthermore, $A \cap [a_1,b_1] = \sn{a_1}$ and $B \cap [a_1,b_1] = \sn {b_1}$.
The pair $\pair{a_1}{b_1}$ is called a \emph{gate} between $A$ and $B$.
The notion of a gate between convex sets is actually defined for a bigger class of compact spaces, called {\em compact median spaces}.
The existence of gates follows from the following fact: Given a family $\Ef$ consisting of closed convex sets with $\bigcap\Ef = \emptyset$, there exist $A, B \in \Ef$ such that $A\cap B = \emptyset$.
For details we refer to the book \cite{Vel}.

\section{Modest semilattices}

In this section we show that $\ceep {K,2}$ is $\sig$-discrete whenever $K$ is a compact totally disconnected semilattice satisfying certain condition.

Namely, we call a semilattice $K$ \emph{modest}, if
it is totally disconnected and for every compact element $p\in K$ the set $p^-$ of all immediate predecessors of $p$ is finite.
Recall that $x$ is an \emph{immediate predecessor} of $p$ if $x < p$ and no $y$ satisfies $x < y < p$.

Below we give two natural examples of modest semilattices.

\begin{prop}
Every compact totally disconnected distributive lattice is a modest semilattice.
\end{prop}

\begin{pf}
Let $\seq{L, \meet, \join, 0,1}$ be such a lattice.
By definition, both operations $\meet$ and $\join$ are continuous, therefore $L$ is a topological semilattice.
Fix $p \in L$ compact and suppose $\sett{x_n}{\ntr} \subs p^-$ is such that $x_n \join x_m = p$ for every $n < m < \nat$.
Let $y$ be an accumulation point of $\sett{x_n}{n>0}$.
From Proposition~\ref{Pnonexplicitsmlts} $y  \ngeq p$, therefore there exists a clopen prime filter $U$ such that $p \in U$ and $y \notin U$.
Find $k < \ell$ such that $x_k, x_\ell \in K \setminus U$.
As $K\setminus U$ is an ideal, $p = x_k \join x_\ell \in K\setminus U$, a contradiction.
This shows that $x \goe p$ and hence $p$ is not a compact element.
\end{pf}

\begin{ex}
Let $T$ be a finitely-branching tree and consider $\al T = T \cup \sn \infty$, where $\infty$ is an additional element satisfying $\infty > t$ for every $t\in T$.
The tree $T$ can be regarded as a locally compact space, where a neighborhood of $t \in T$ is of the form $(s,t]$, where $s < t$.
Define the topology on $\al T$ so that it becomes the one-point compactification of $T$.
Define an operation $\meet$ on $\al T$ as follows:
let $s \meet t = \max(s,t)$ whenever $s$ and $t$ are comparable, and let $s\meet t = \infty$ otherwise.
It is easy to check that this is a continuous semilattice operation.
The property that $T$ is finitely-branching is equivalent to the fact that $\al T$ is a modest semilattice.
\end{ex}

\begin{tw}
Let $K$ be a modest semilattice. Then $\ceep {K,2}$ is $\sig$-discrete.
\end{tw}

\begin{pf}
Throughout this proof, we shall consider $\dsl K$ as the set of all positive compact elements, plus an artificial element $\infty\notin K$ which corresponds to the empty clopen filter.
So the meet operation on $\dsl K$ is actually the supremum in $K$.
Given $p,q\in \dsl K$ we shall denote by $p\cdot q$ the meet of $p$ and $q$ in $\dsl K$, which equals either $\sup\dn pq$ in $K$ or $\infty$ in case $\dn pq$ is not bounded from above in $K$.

Fix $f \in \cee K$.
We say that $f$ has a \emph{jump} at $p\in \dsl K$ if $p \ne \infty$ and $f(p)\ne f(q)$ for every $q\in p^-$.

\begin{claim}
Assume $f\in \cee {K,2}$. If $f$ is not constant then $f$ has a jump at some $p\in \dsl K$.
\end{claim}

\begin{pf}
Let $p$ be a minimal element of $K$ such that $f(p) \ne f(0)$.
Then $p\in \dsl K$ by Proposition~\ref{Pnonexplicitsmlts}(4) and clearly $f$ has a jump at $p$.
\end{pf}

We shall say that $f\in \cee K$ has a \emph{relative jump} at $p\in \dsl K$ \emph{with respect to} $q\in\dsl K \cup \sn 0$ if $q < p$ and $f(x)\ne f(p)$ for every $x\in p^- \cap [q)$.
Note that a relative jump with respect to $0$ is just a jump.

Fix $f\in \cee{K,2}$. Let $L_0(f)\subs \dsl K$ be the subsemilattice generated by all $p\in \dsl K$ such that $f$ has a jump at $p$.

By induction, we define $L_n(f)$ to be the subsemilattice of $\dsl K$ generated by $L_{n-1}(f)$ together with
all $p\in \dsl K$ such that $f$ has a relative jump at $p$ with respect to some $q\in L_{n-1}(f)$.

\begin{claim}\label{Cwejiq}
For every $f\in \cee {K,2}$ there exists $\ntr$ such that $L_n(f) = L_{n+1}(f)$.
\end{claim}

\begin{pf}
Let $\map \phi KF$ be a continuous epimorphism onto a finite semilattice such that $f$ is constant on each fiber of $\phi$.
Note that $\dsl F \subs \dsl K$, after the obvious identification (the pre-image of a compact element is compact).
We shall prove by induction that $L_n(f) \subs \dsl F$ for every $\ntr$.

Notice that a relative jump with respect to $0$ is just a jump.
Thus, set $L_{-1}(f) = \sn \infty$, in order to start the induction.

Assume now that $L_{n-1}(f) \subs \dsl F$ and fix $p\in L_n(f)$ such that $f$ has a relative jump at $p$ with respect to $q\in L_{n-1}(f)$.
Find $t\in F$ such that $p\in \phi^{-1}(t)$.
Let $u = \min \phi^{-1}(t) \in \dsl F$.
We need to show that $p = u$.
For suppose that $u < p$.
Find $r\in p^-\cap [u)$.
Notice that $q \loe u$, because $[q)$ is the union of some fibers of $\phi$ and $p$, $u$ are in the same fiber.
It follows that $f$ does not have a relative jump at $p$ with respect to $q$, because $f(p)=f(r)$.
This is a contradiction, which shows that $L_n(f) \subs \dsl F$.
\end{pf}

Let $\sig$ denote the set of all finite sequences of natural numbers.
Given $f\in \cee {K,2}$, let $n(f)$ be the minimal $\ntr$ such that $L_{n+1}(f) = L_n(f)$.
By Claim~\ref{Cwejiq}, this is well defined.
Furthermore, let $s(f)(i) = |L_i(f)|$ for $i \loe n(f)$.
So $s(f) \in \sig$.

This provides a decomposition of $\cee {K,2}$ into countably many pieces, indexed by $\sig$.
Namely, we shall prove that given $s\in\sig$, the set
$$C_s = \setof{f\in \cee {K,2}}{s(f) = s}$$
is discrete with respect to the pointwise convergence topology.
In order to avoid too many indices, we set $L(f) = L_{n(f)}(f)$.

The following claim is crucial.

\begin{claim}\label{Curoqe}
Let $f\in\cee{K,2}$ and let $\map \phi KS$ be the continuous semilattice homomorphism induced by the inclusion $L(f) \subs \dsl K$ (i.e. $S = \dsl {L(f)}$ and $\phi$ is the dual map to this inclusion).
Then $f$ is constant on each fiber of $\phi$.
\end{claim}

\begin{pf}
Note that for each $t\in S$ the minimal element of $\phi^{-1}(t)$ belongs to $L(f)$.

Suppose that $f$ is not constant on $G = \phi^{-1}(t)$ and let $q=\min G$.
Then $f$ has a jump at some compact element $p\in G$, $p>q$.
This jump is relative with respect to $q$.
But $q\in L_n(f)$ for some $n\loe n(f)$ so $p\in L_{n+1}(f) \subs L(f)$.
We conclude that $q$ is the minimal element of some fiber of $\phi$, a contradiction.
\end{pf}

Fix $f\in \cee {K,2}$ and define
$$M_f = L(f) \cup \bigcup_{p\in L(f)}p^-.$$
Since $K$ is modest, $M_f$ is finite.
Assume $g\in \cee {K,2}$ is such that $f\rest M_f = g\rest M_f$ and $s(f) = s(g)$.
We claim that $f=g$.

For this aim, we shall prove by induction that $L_n(g) = L_n(f)$ for every $n\loe n(f)=n(g)$.
This is trivial for $n=-1$.
In fact, this also holds for $n=0$, because given $p\in L_0(f)$ for which $f$ has a jump, $g$ has the same values as $f$ in the finite set $\sn p \cup p^-$.
It follows that $g$ has a jump at $p$.
Hence $L_0(g)\subs L_0(f)$ and consequently $L_0(g) = L_0(f)$ because these sets have the same finite cardinality.

Now assume that $L_{n-1}(g) = L_{n-1}(f)$ and fix $p\in L_n(f)$ such that $f$ has a relative jump at $p$ with respect to some $q\in L_{n-1}(f)$.
Then $q\in L_{n-1}(g)$ and $g$ has exactly the same values as $f$ on the finite set $\sn p \cup (p^-\cap [q))$.
It follows that $g$ has a relative jump at $p$ with respect to $q$, i.e. $p\in L_n(g)$.
This shows that all generators of $L_n(f)$ are contained in $L_n(g)$.
Thus $L_n(f) \subs L_n(g)$ and consequently $L_n(f) = L_n(g)$, again because of the same finite cardinality.

Finally, by Claim~\ref{Curoqe}, we conclude that $g = f$, because $M_f$ intersects each fiber of the quotient map induced by $L(f)$.
\end{pf}

Actually, the proof above provides a recipe for the minimal semilattice quotient for a given function into $\dn 01$.

\section{Finite-dimensional compact lattices}

In this section we prove that $C(K)$ has the {\bf SLD} whenever $K$ is a finite-dimensional compact lattice.
Recall that a lattice $K$ is \emph{finite-dimensional} if $K \subs L_1 \times \dots \times L_n$, where each $L_i$ is a compact totally ordered space (called briefly a \emph{compact line}).
Note that every compact totally ordered space is a continuous image of a 0-dimensional one, therefore we can restrict attention to 0-dimensional lattices.

A result from \cite{BurKuTod} says that $C(K)$ has a pointwise Kadets renorming (which implies {\bf SLD}) whenever $K$ is a product of compact lines.
Up to now, it is not known whether the same result holds for $K$ being a closed sublattice of a finite product of compact lines.
This question remains open, however we prove a weaker result concerning the {\bf SLD} property.

From now on $L_i$, $1\leq i\leq n$, will stand for totally ordered spaces that are zero--dimensional compact for their order topology.  Let $\leq$ be the {\em product order} in $\prod_{i=1}^{n}L_{i}$, i.e. $\left( a_{i}\right)_{i=1}^{n} \geq  \left( b_{i}\right)_{i=1}^{n}$ whenever $a_{i} \geq b_{i}$ for any $i\leq n$; this space endowed with it and the product topology is a zero--dimensional compact distributive lattice.

Until the the end of the proof of Proposition~\ref{dim22020} we assume that $K$ is a fixed sublattice of $\prod_{i=1}^{n}L_{i}$, endowed with a topology ${\mathcal T}$ for which it is a compact distributive lattice that, by compactness, must be the subspace topology of $\prod_{i=1}^{n}L_{i}$.

Given $x=\left( x_{i}\right)_{i=1}^{n}\in \prod_{i=1}^{n}L_{i}$ let $\pi_{i} (x)=x_i$, and $\rho_{i} (x)=\left( x_{j}\right)_{j\neq i}$, $1\leq i\leq n$.  If $1\leq i\leq n$, and $p\in \Pri   \left( L_{i} \right)$ we write $\tilde{p} := \left( \prod_{j<i}L_{j}\times p\times \prod_{j>i} L_{j}\right) \cap K $.  It is clear that $\tilde{p}  \in \Pri   \left( K\right)$.

\begin{tw}\label{dim21937}
$C(K)$ has the property {\bf SLD} for the pointwise topology.
\end{tw}

It will be deduced from some results.

\begin{lm}\label{dimn1135}
For $1\leq i\leq n$ let $p$, $q\in \Pri   \left( L_{i} \right)$ with $p\sups q$, if $\pair ab$ is a gate between $\pair  {\tilde{p}}{\tilde{q}}$ we have $\pi_{i}(a)< \pi_{i}(b)$; and for any $c\in K$ we have $c\notin \tilde{q}\cup (K\setminus \tilde{p})$ whenever $\pi_{i}(a) <\pi_{i}(c)< \pi_{i}(b)$.  In particular $K\setminus \tilde{p}=\left\{ x\in K:\ \pi_{i}(x) \leq  \pi_{i}(a) \right\}$ and $\tilde{q}=\left\{ x\in K:\ \pi_{i}(x) \geq  \pi_{i}(b) \right\}$.  Moreover for any $c\in K$ we have
\begin{eqnarray}
&&\hbox{if}\  \pi_{i}(c)\geq \pi_{i}(b)\ \hbox{\rm and}\ \pi_{j}(a)<\pi_{j}(b) \ \hbox{then}\ \pi_{j}(c)\geq \pi_{j}(b)\ 1\leq j\leq n  .\label{dim2007}
\end{eqnarray}
\end{lm}
\begin{pf}
For the first part, from $\pi_{i}((a\vee c)\wedge b) =\pi_{i}(c)$ we get $ (a\vee c)\wedge b\in [a,b] \setminus \{ a,b\}$.  Then, since $\pair ab$ is a gate we have $(a\vee c)\wedge b\notin \tilde{q}\cup (K\setminus \tilde{p})$.   But $\tilde{q}$ is a filter so $c\notin \tilde{q}$.  Similar arguments show that $c\notin K\setminus \tilde{p}$.

To show (\ref{dim2007}) observe that since $\pair ab$ is a gate we have $a<b$.  Suppose that $c$ and $j$ witness that (\ref{dim2007}) is false.  We have $c\in \tilde{q}$ and $a\vee (b\wedge c)\in \tilde{q}$.  Since $\pair ab$ is a gate we have $a\vee (b\wedge c) =b$, so $\pi_{j}( b ) =\pi_{j}(a)\vee (\pi_{j}(b)\wedge \pi_{j}(c))) = \pi_{j}(a)\vee   \pi_{j}(c)$ which contradicts that $\pi_{j}(a)<\pi_{j}(b)$ and $\pi_{j}(c)<\pi_{j}(b)$.
\end{pf}

\begin{lm}\label{dim21816}
Let $a$, $b\in K$ be \st\ $\left\{ x\in K\cap [a,b]:\ \pi_{i}(a) <\pi_{i}(x)<\pi_{i}(b)\right\}=\emptyset$, then $\left\{ x\in K:\ \pi_{i}(a) <\pi_{i}(x)<\pi_{i}(b)\right\}=\emptyset$.
\end{lm}

\begin{pf}
Since $(a\vee x)\wedge b\in [a,b]$ for any $x\in K$ and $\pi_{i}((a\vee x)\wedge b)=\pi_{i}(x)$ the assertion follows.
\end{pf}

If $1\leq i\leq n$, $p$, $q\in \Pri   \left( L_{i} \right)$
with $p \sups q$, $\pair ab$ is a gate between
$\pair {\tilde{p}}{\tilde{q}}$ and $a$, $b$ belong to a subset $S\subs K$ we will say that $f$ $(i,\eps)${\em --leaps} (resp. {\em jumps}) at $\pair {\tilde{p}}{\tilde{q}}$, or at $\pair ab$ or within $S$ whenever it $ \eps$--leaps (resp. jumps) at $\pair {\tilde{p}}{\tilde{q}}$.

\begin{df}
\label{dim21330}
{\em Given $f\in C(K)$, $1\leq i\leq n$, and $\eps >0$ we will say that $f$ {\em $(i,\eps)$--leaps $m$ times} if there exists $\left\{ \left< p_{i,j},q_{i,j}\right> \right\}_{j=1}^{m}$, $p_{i,j}$, $q_{i,j}\in \Pri   \left( L_{i} \right)$ \st\ $\tilde{p}_{i,1}\sups \tilde{q}_{i,1} \supsetneqq \tilde{p}_{i,2}\sups \tilde{q}_{i,2} \supsetneqq  \ldots \supsetneqq \tilde{p}_{i,m}\sups \tilde{q}_{i,m}$; and
\begin{itemize}
\item[(i)]  $f$ $\eps$--leaps at each $\left( \tilde{p}_{i,j},\tilde{q}_{i,j}\right)$, $1\leq j \leq m$;
\item[(ii)] if $f$ $\eps$--jumps at $\tilde{p}$ for some $p\in \Pri   \left( L_{i}\right)$, $1\leq i\leq n$, then there exists $j$, $1\leq j \leq m$, \st\ $p=p_{i,j}=q_{i,j}$.
\end{itemize}
}
\end{df}
From Lemma~\ref{dim21130} below it follows that this definition makes sense.
\begin{uwgi}\label{dim22009}
According to Lemma~\ref{dimn1135} given $1\leq i\leq n$ and $p_k$, $q_{k}\in \Pri   \left( L_{i} \right)$, $k=1$, $2$, \st\ $\tilde{p}_{1}\sups \tilde{q}_{1} \supsetneqq \tilde{p}_{2} \sups  \tilde{q}_{2} $, and $\pair {a_{j}}{b_{j}}$ gates between $\pair {\tilde{p}_{j}}{\tilde{q}_{j}}$, $j=1$, $2$, we have $\pi_{i}\left( a_{1} \right) <\pi_{i}\left( b_{1} \right) \leq  \pi_{i}\left( a_{2} \right)$.
\end{uwgi}
\begin{lm}\label{dim21130}
Given $f\in C(K)$ and $\eps >0$ the set of all $p\in \Pri   \left( L_{i} \right)$ \st\ $f$ $(i,\eps)$--jumps at $\tilde{p}$ is finite and so is the set of all $m\in \N$ \st\ $f$ $(i,\eps)$--leaps $m$--times.
\end{lm}
\begin{pf}
By compactness for any $1\leq k\leq n$ it is possible to choose $\left\{ \alpha_{k,j}\right\}_{j=0}^{\ell_k}$, $\alpha_{k,j}\in L_k$, \st\ the $\left[ \alpha_{k,j-1} ,\alpha_{k,j} \right]$'s cover $L_k $ and the oscillation of $f$ on each $K\cap \prod_{k=1}^{n} \left[  \alpha_{k,j -1}, \alpha_{k,j }\right]$ is strictly less than $\eps$; those sets cover $K$. Let $p_m$, $q_{m}\in \Pri   \left( L_{i} \right)$, $m=1$, $2$, \st\ $f$ $(i,\eps)$--leaps at each $\left( \tilde{p}_{m},\tilde{q}_{m}\right)$, with $\tilde{p}_{1}\supset\tilde{q}_{1} \supsetneqq \tilde{p}_{2} \supset \tilde{q}_{2} $, and this is witnessed by the gates $\pair {a_{m}}{b_{m}}$, $m=1$, $2$.  We have that $a_{m}$ and $b_{m}$ cannot be in the same  $\prod_{k=1}^{n} \left[  \alpha_{k,j -1}, \alpha_{k,j }\right]$.  To finish the proof it is enough to apply Remark~\ref{dim22009} and (\ref{dim2007}).
\end{pf}
Fix  $f \in C(K)$ and $\eps >0$ until the end of the proof of Proposition~\ref{dim22020}.  Given $1\leq i\leq n$ let
$m_{i}=m_{i}(f)$ the maximum number of times that $f$ $(i,\eps)$--leaps, so let $\left\{ \left< p_{i,j},q_{i,j}\right> \right\}_{j=1}^{m_i}$ satisfying (i) and (ii) of Definition~\ref{dim21330}.  Let $C_{k}:=\prod_{i=1}^{n}\left[  a_{k,i}, b_{k,i}\right]\cap K$, $1\leq k\leq \ell$, $a_{k,i}$, $b_{k,i}\in L_i$, any covering of $K$ made up by sets \st\ each $C_{k}$ is either included or disjoint of every one of the sets $K\setminus \tilde{p}_{i,j}$, $\tilde{q}_{i,j}$, $1\leq j\leq m_i$, $1\leq i\leq n$.  From now on we will fix $k$ and, for simplicity, we will write $C_{k} =\prod_{i=1}^{n}\left[  a_{i}, b_{i}\right]\cap K$.  Let the sequence $p(n)$ defined inductively by $p(1)=3$, $p(n+1)=8p(n)+6$, $n\in \N$.
\begin{prop}\label{dim22020}
${\displaystyle \osc \left( f,  C_{k} \right) \leq p(n)\eps .}$
\end{prop}
It will be a consequence of some lemmata.
\begin{lm}\label{dim20927}
Given $u$, $v \in K$, $1\leq i\leq n$, $U$, $V$ open sets in $L_i$ with $|f(u)-f(v)| >\eps$, $\rho_{i}(u)=\rho_{i}(v)$, $\pi_{i}(u) <\pi_{i}(v)$ and $\pi_{i}(u)\in U$, $\pi_{i}(v)\in V$, there exist $a$, $b \in [u,v]\cap K$ with $\pi_{i}( a)\in U$, $\pi_{i}(b)\in V$, $ \pi_{i}(a) <\pi_{i}(b) $, \st\ $\pair ab$ witnesses an $(i,\eps)$--leap of $f$.
\end{lm}
\begin{pf}
Since $S:=\left\{ x\in K:\ \rho_{i}(x)=\rho_{i}(u)\right\}$ is a $0$--dimensional totally ordered compact space, the right (resp. left) isolated points are dense in it.  From the continuity of $f$ we get two open sets $U$, $V$ in $L_i$ and a gate $\pair ab$ that witnesses an $\eps$--leap of the restriction of $f$ to $S$.  Lemma~\ref{dim21816} shows that $\pair ab$ is a gate in $K$ too.
\end{pf}
\begin{lm}\label{dim21832}
Let $a$, $b \in C_k$,  if $M:=\max f\left( [a,b] \cap K\right)$ and $m:=\min  f\left( [a,b]  \cap K\right)$ then $[M,m]\setminus f\left([a,b] \cap K\right)$ contains no interval of length bigger that $\eps$.
\end{lm}
\begin{pf}
By contradiction let $r$, $s\in \R$ \st\ $r<s$, $s-r>\eps$ and $[r,s] \cap f(  [a ,b ] \cap K) =\{ r,s\}$; We must have that $D:=\{ i:\ 1\leq i\leq n,\ \pi_{i}(a) \neq \pi_{i}(b)\}$ is non empty.  We argue by induction over $\ell:=\hbox{card~}D$.  Let $v$, $w\in [a ,b ] \cap K$, \st\ $f(v)=s$, $f(w)=r$.  When $\ell =1$ the set $[v,w]$ is totally ordered.  If, for instance, $v<w$, set $v_{0}:=\sup \{ y \in [v ,w ] \cap K: \ f(y )\geq s  \} $ and $w_{0}:=\inf \{ y\in \left[ v_{0},w \right] \cap K:\  f(y )\leq r \}$.  Thus $f\left( v_{0} \right) -f\left( w_{0} \right)\geq s -r >\eps$ and $\left[  v_{0} ,  w_{0} \right] \cap K   =\left\{  v_{0} , w_{0} \right\}$.  Lemma~\ref{dim21816} shows that $f$ $(i,\eps)$--jumps at $\left(  v_{0} , w_{0} \right)$.  A contradiction.
\arre
Assume that the assertion holds for $k$, $1\leq k<n$, and $\ell =k+1$.  Fix $i$, $1\leq i\leq n$, \st\ $\pi_{i} (v)\neq  \pi_{i} (w)$, and, for instance, $\pi_{i} (v)<  \pi_{i} (w)$.  Now set $v_{0}$ \st\ $ \pi_{i} \left( v_{0}\right)=\sup \{ \pi_{i} (y) :\ y \in [v ,w ] \cap K,  \ f(y )\geq s  \} $ and $w_{0}$ \st\ $ \pi_{i} \left( w_{0}\right)=\inf \{ \pi_{i} (y) :\ y \in \left[ v_{0} ,w \right] \cap K,  \ f(y )\leq r  \} $.  From the choice of $\left( v_{0},w_{0}\right)$ and Lemma~\ref{dim21816} if follows that $\left\{ x\in K:\ \pi_{i} \left( x\right) \leq \pi_{i} \left( v_{0}\right)\right\}$ and $\left\{ x\in K:\ \pi_{i} \left( x\right) \geq \pi_{i} \left( w_{0}\right)\right\}$ are a proper clopen ideal and a filter respectively, whose union is $K$.  On the other hand the sets
$$
\left\{ x\in [v_{0},w_{0}] \cap K:\ \pi_{i} \left( x\right)  = \pi_{i} \left(  v_{0}\right) \right\} \ \hbox{and}\ \left\{ x\in [v_{0},w_{0}] \cap K:\ \pi_{i} \left( x\right)  = \pi_{i} \left(  w_{0}\right) \right\}
$$
are nonempty disjoint closed convex subsets of $K$, then there exists a gate $\pair {\alpha}{\beta}$ between them, it is clear that it must be a gate between the clopen ideal and filter above and $\pi_{i} \left(  \alpha\right) =\pi_{i} \left(  v_{0}\right)$, $\pi_{i} \left(  \beta\right) =\pi_{i} \left(  w_{0}\right)$.  We get that the cardinal of the set $\left\{ j:\ 1\leq j\leq n,\ \pi_{j}(\alpha ) \neq \pi_{j}\left(  v_{0}\right)\right\}$ is not bigger than $k$ so, according to the inductive hypothesis applied to $\left[ \alpha ,v_{0}\right]$ we have $f(\alpha )\geq s$; the same argument shows that $f(\beta )\leq r$.  Then we have found an $\eps$--jump which is a contradiction.
\end{pf}
\noindent
{\em Proof of Proposition~\ref{dim22020}.}  If $a$, $b\in C_k$ and $\{i:\ 1\leq i\leq n, \pi_{i}(a)\neq \pi_{i}(b)\}$ has cardinal $k$, we will show by induction on $k$ that $|f(a)-f(b)| <p(k)\eps$.  Assume that $k=1$.  By contradiction let $1\leq i\leq n$ with $\rho_{i}(a)=\rho_{i}(b)$ and $|f(a)-f(b)| > 3\eps$.  According to Lemma~\ref{dim21832} there exists $\xi\in [a,b]\cap K$ \st\ $|(1/2)(f(a)+f(b))-f(\xi)| <\eps$.  Then $| f(a) -f(\xi)| >\eps$ and $| f(b )-f(\xi)| >\eps$.  From Lemma~\ref{dim20927} we get two different $(i,\eps)$--leaps.  A contradiction.\arre
Assume that the assertion holds for $k$, $1\leq k<n$, and $\ell =k+1$.  By contradiction suppose that $|f(a)-f(b)| > p(k+1)\eps$ and $\ell=k+1$.  Since $p(k+1)\eps <|f(a)-f(b)|\leq |f(a)-f(a\vee b) |+ |f (a\vee b)-f(b)|$ we may assume that $a<b$ and $|f(a)-f(b)| >(1/2)p(k+1)\eps$.  From Lemma~\ref{dim21832} we get $c\in [a,b]\cap K$ \st\ $|f(a)-f(c)|>(1/2)((1/2)p(k+1)-1)\eps$ and $|f(c)-f(b)|>(1/2)((1/2)p(k+1)-1)\eps$.  We claim that this implies that $f$ $(i,\eps)$--leaps in $[a,c]$ and in $[c,b]$ which contradicts the maximality of $m_i$.  Let us show it in $[a,c]$.  Since $L_i$ is $0$--dimensional, the continuity of $f$ allows us to assume that $\pi_{i}(a)$ (resp. $\pi_{i}(c)$) is right (resp. left) isolated.  Then $\left\{ x\in [a,b]:\ \pi_{i}(x)= \pi_{i}(a)\right\}$ and $\left\{ x\in [a,b]:\ \pi_{i}(x)= \pi_{i}(c)\right\}$ are nonempty clopen convex sets; let $\pair uv$ a gate between them,  that should be a gate between the clopen ideal and filter $\left\{ x\in [a,b]:\ \pi_{i}(x)\leq \pi_{i}(a)\right\}$ and $\left\{ x\in [a,c]:\ \pi_{i}(x)\geq \pi_{i}(c)\right\}$ too.  It is clear that $\pi_{i} \left(  u\right) =\pi_{i} \left(  a\right)$, $\pi_{i} \left(  v\right) =\pi_{i} \left(  c\right)$ and the cardinal of the sets $\left\{ j:\ 1\leq j\leq n,\ \pi_{j}(u ) \neq \pi_{j}\left(  a\right)\right\}$ and $\left\{ j:\ 1\leq j\leq n,\ \pi_{j}(v ) \neq \pi_{j}\left(  c\right)\right\}$ are not bigger than $k$ so, according to the inductive hypothesis $|f(a)-f(u)|\leq p(k)\eps$ and $|f(v)-f(c)|\leq p(k)\eps$.  Then $|f(u)-f(v)|> ((1/2)((1/2)p(k+1)-1)-2p(k))\eps\geq \eps$ so $\pair uv$ witnesses an $(i,\eps)$--leap.\hfill $\square$

{\em Proof of Theorem~\ref{dim21937}}  Given $\eps >0$ we can write $C(K)=\bigcup_{p=1}^{\infty} C_p$ in such a way that for any $p\in \N$ there exist $m_i$, $1\leq i\leq n$ \st\ $m_{i}(f)=m_{i}$, $1\leq i\leq n$, for al $f\in C_p$.  For a $f\in C_p$ these leaps are witnessed by filters and gates $\pair {a_{i,j}}{b_{i,j}}$, $1\leq j\leq m_i$, $1\leq i\leq n$, with $\left| f\left( a_{i,j}\right) -f\left( b_{i,j}\right) \right| >\eps$.  According to Proposition~\ref{dim22020} if
$$
g\in \left\{ h\in C_{p}:\ \left| h\left( a_{i,j}\right) -h\left( b_{i,j}\right) \right| >\eps,\ 1\leq j\leq m_i, 1\leq i\leq n \right\}
$$
then $\osc \left( g,  K\cap \prod_{i=1}^{n}\left[  a_{i}, b_{i}\right]  \right) \leq p(n)\eps$.  To finish apply \cite[Theorem~4.(iii)]{JFtopol}.
\hfill $\square$

%

{\sc W. Kubi\'s, Institute of Mathematics, Czech Academy of Sciences, \v{Z}itn\'a 25, 115 67 Praha 1,
CZECH REPUBLIC} \textit{and} {\sc Institute of Mathematics, Jan Kochanowski University in Kielce, POLAND}

\emph{E-mail address:} \texttt{kubis@math.cas.cz}\medskip

{\sc A. Molt\'o, Departamento de An\'alisis Matem\'atico, Facultad
de Ma\-te\-m\'aticas, Universidad de Valencia, Dr. Moliner 50, 46100
Burjassot (Valencia), SPAIN}

\emph{E-mail address:} \texttt{anibal.molto@uv.es}\medskip

{\sc S. Troyanski, Departamento de Matem\'aticas, Universidad de Murcia, Campus de Espinardo, 30100
Espinardo(Murcia), SPAIN}

\emph{E-mail address:} \texttt{stroya@um.es}

\end{document}